 \documentclass{amsart}

  \let\oldchapter=\chapter
\def\resetpage{\setcounter{page}{1}}
\def\chapter{\expandafter\resetpage\oldchapter}  
 
\usepackage{mathdots}

\usepackage{bm}
\usepackage{bigints}
\usepackage{epsfig}
\usepackage{epstopdf}

\usepackage{srcltx,amssymb}

\usepackage{graphicx}

    \usepackage{amsthm}
      \usepackage{amsmath}
    \usepackage{amscd}
    \usepackage[title]{appendix}
    \usepackage{slashed}
    \numberwithin{equation}{section}

   \usepackage{appendix}

    \theoremstyle{definition}

   \numberwithin{theorem}{section}

\begin{document}

        \title{Cosmology from a non-physical standpoint: an algebraic analysis}
    \author{Fred Greensite}

  \begin{abstract} We present a non-physical interpretation of the Cosmological Constant based on a particular algebraic analysis. This also introduces some novel algebraic structures, such as ``unital norms", ``uncurling metrics", and ``partial wedge products".  \end{abstract}
  
  \maketitle
  \markright{Cosmology from a non-physical standpoint}

\section{Introduction}

Based on independent astronomical observations from the 1990s, it is generally believed that the expansion of the universe is accelerating  \cite{riess:1998,perlmutter:1999}.  Within the formalism of general relativity, this is accommodated by addition of the term  $\Lambda g_{\mu\nu}$ to the left-hand-side of the original Einstein equations, where  $g_{\mu\nu}$ is the metric tensor of the pseudo-Riemannian spacetime manifold, and the Cosmological Constant $\Lambda$ has a small positive value with respect to usually applied systems of physical units.
The resulting Einstein equations are, 
  \begin{equation}\label{7/15/18.1}  R_{\mu\nu} -\frac{1}{2}R \,g_{\mu\nu}+ \Lambda g_{\mu\nu} = \left(\frac{8\pi G}{c^4}\right)T_{\mu\nu}, \end{equation} where $R_{\mu\nu}$ and $R$
are respectively the Ricci curvature tensor and scalar curvature, $G$ is Newton's gravitational constant, $c$ is the speed of light, and $T_{\mu\nu}$ is the stress-energy tensor.  Examining (\ref{7/15/18.1}) in the context of the usual spacetime $4$-vector $x^\mu = (x^0,x^1,x^2,x^3)$ whose components have physical units of length (so that the metric tensor is devoid of physical units), the curvature expressions on the left-hand-side of (\ref{7/15/18.1}) indicate that $\Lambda$ has units of length$^{-2}$ ($x^0\equiv ct$ with $t$ representing ``time").
 To emphasize its role as a source of dynamics, $\Lambda g_{\mu\nu}$ is sometimes instead subtracted from the right-hand-side of the original Einstein equations.     But a third equivalent presentation is removal of $\Lambda$ as a coefficient in the new additional term and instead subsuming it within the curvature expressions, initiated by simply dividing both sides of (\ref{7/15/18.1}) by $\Lambda$.  For example, it can then be incorporated into  the default coordinate vector field $\frac{\partial}{\partial x^\mu}$ underlying the curvature expressions, so that coordinate vector field becomes $\frac{\partial}{\partial(\sqrt{\Lambda}x^\mu)}$.  This third format highlights the fact that a non-vanishing Cosmological Constant renders the Einstein equations unitless.
 
  In other words, for the Einstein equations augmented by the Cosmological Constant, the unitless metric tensor $g_{\mu\nu}$ as an entity varying on the unitless domain with points $\sqrt{\Lambda}x^\mu$, is dependent on the unitless source $\left(\frac{8\pi G}{c^4\Lambda}\right)T_{\mu\nu}$.  Of course, the point is that not only do no physical units appear in the equations,  but $g_{\mu\nu}$,   $\sqrt{\Lambda}x^\mu$, and $\left(\frac{8\pi G}{c^4\Lambda}\right)T_{\mu\nu}$, are not functions of the choice of a system of physical units.

This elimination of all physical units might be viewed as an exoneration of the Pythagorean conception of a purely mathematical essence comprising nature (at least insofar as gravitation theory is concerned).  There seems to be an inexorable historical progression in that direction.  \begin{itemize}\item The experimental analyses of {\it Galileo} made it possible to talk scientifically about three physical units: those of mass, length, and time.  \item With subsequent recognition of the Lorentz covariance of {\it Maxwell}'s equations, those three basic units reduce to two via the universal constant $c$, e.g., only mass (i.e., energy) and length units remain.\item In the context of the universal constant $G$ first introduced in {\it Newton}'s law of gravitation, those two basic units reduce to one, e.g., only the length unit remains (characterizing the geometric unit system of general relativity). \item   Introduction of the (presumably) universal constant $\Lambda$, arising from {\it Einstein}'s recognition of an augmented equation set satisfying the general covariance principle, then reduces the number of physical units from one to none. \end{itemize}

Together, $c, G,\Lambda$ scrub the Einstein equations clean of physical units.  So, the question arises:  \begin{quote} If physical units need not appear in the basic equations of general relativity, why should they appear in the equations of special relativity? \end{quote}

    To add perspective, imagine that we instead make the {\it prior} assumption that physical units are ultimately not fundamental in physics (perhaps something Plato might have liked).  Once the Einstein equations are derived in the usual physical way (where the means of satisfying Einstein's general covariance principle is by equating the stress-energy tensor to an expression that is required to be linear in the Riemann tensor and dependent solely on the Riemann tensor and metric tensor),  the {\it prediction} will be that $\Lambda$ is nonzero, so that physical units are eliminated.  Given that general relativity is already based on the (equivalence) principle indicating that gravitational effects are geometrical, the only thing missing from such a non-physical program would be a derivation of the Lorentzian geometry of the spacetime manifold tangent spaces in a manner in which physical units play no part.  
    
    To add credence to the prediction that the Einstein equations should be devoid of physical units, it is of interest to present a program wherein physical units {\it cannot} have a fundamental role in special relativity.  This would be effected by deriving the tangent space at each point of the spacetime manifold as a self-contained algebra (mathematically, a tangent space at a point of the spacetime manifold is intrinsically indistinguishable from the space of events in special relativity with respect to a particular origin).  Physically, the members of the tangent space must carry the same physical unit(s).  But the product of two algebra members would then be attached to the square of the physical unit(s), and so could not map back into the tangent space.  Thus, if a spacetime manifold tangent space is a self-contained algebra, physical units must ultimately disappear.
         
     So, suppose we actually take this line of thought seriously. Can we derive the Minkowski spacetime metric in this manner from only the structural components of an inertial observer's physical space as given to us by Euclid, and can anything be learned from that?
    We intend to answer those questions in the positive, by presenting such a derivation that provides a broader context for considering the emergence/significance of spacetime geometry.  More fancifully, the results could also be considered in the context of currently speculative notions such as ``histories exist, rather than time", or ``we are living in a simulation". 
    
\section{Four-dimensional spacetime algebras} 
\subsection{The unital case} According to the penultimate paragraph of the prior section we want to identify a four-dimensional spacetime algebra to accomplish our objective of a unitless special relativity.  We begin with the given algebra of real numbers $\mathbb{R}$ and imposition of the property of linearity whenever the opportunity arises ($\mathbb{R}$ itself is  the very embodiment of linearity through its product and sum operations).  An ``observer" notes that locally his surrounding ``space" seems well modeled by $\mathbb{R}^3$ (taken to be a linear space, i.e., a vector space), and  the fact that the product and sum of $\mathbb{R}$ with $\mathbb{R}^3$ also seem to have descriptive power.  That is, the tensor product $\mathbb{R}\otimes\mathbb{R}^3$, a linear space, models the feature of ``scale" through the component-wise product of scalars with points of $\mathbb{R}^3$.  The direct sum $\mathbb{R}\oplus\mathbb{R}^3$, a linear space, allows indexing of serial configurations (trajectories) in his physical space.  He also notices that as a linear space,  $\mathbb{R}^3\approx\mathbb{R}\otimes\mathbb{R}^3$ seems to conform to Euclidean geometry (so it is an inner product space).  In what follows, we will take $\mathbb{R}^{m,n}$ to be the vector space $\mathbb{R}^{m+n}$ associated with a symmetric bilinear form having signature $(m,n)$.  Thus, the observer finds his physical space to be $\mathbb{R}^{3,0}$.    Given the known geometry of his physical space, the observer is interested in geometry of the augmented space $\mathbb{R}\oplus\mathbb{R}^{3,0}$  arising from the need to describe trajectories.
 
 Since our program can be placed within the context of Geometric (Clifford) Algebras, we will initially generalize the notation by writing the observer's physical space $\mathbb{R}^{3,0}$ as $\mathbb{R}^{m,n}$.   
 Following the variant nomenclature commonly used with the latter algebras, we refer to the bilinear form associated with $\mathbb{R}^{m,n}$ as an inner product.

 We first note that there are three bilinear products underlying the linear space $\mathbb{R}^{m,n}$, these being the field product $\mathbb{R}\times \mathbb{R}\rightarrow \mathbb{R}$, component-wise multiplication by a scalar, $\mathbb{R}\times \mathbb{R}^{m+n}\rightarrow \mathbb{R}^{m+n}$, and the inner product $\mathbb{R}^{m+n}\times \mathbb{R}^{m+n}\rightarrow \mathbb{R}$.  Linear extrapolation  implies a bilinear product ``$\,\bullet\,$" making $\mathbb{R}\oplus\mathbb{R}^{m+n}$ an algebra such that for
$\alpha,\beta\in \mathbb{R},\,\, {\sf a},{\sf b}\in \mathbb{R}^{m+n}$,  \begin{eqnarray}\nonumber  (\alpha+ {\sf a})\bullet (\beta + {\sf b}) &=& \alpha\bullet\beta + \alpha\bullet {\sf b} + {\sf a}\bullet\beta + {\sf a}\bullet {\sf b} \\&=&
    [\alpha \beta + {\sf a} \cdot {\sf b}] + [ \beta {\sf a} + \alpha {\sf b}] \in \mathbb{R}\oplus \mathbb{R}^{m+n}, \label{12/24/18.1}\end{eqnarray} where  ``\,$\cdot\,$" is the inner product on $\mathbb{R}^{m,n}$.  That is, the three products inherent in $\mathbb{R}^{m,n}$ imply a groupoid structure on $\mathbb{R}\cup \mathbb{R}^{m+n}$ with a product $\bullet$ such that $\alpha\bullet\beta\equiv \alpha \beta$, $\alpha \bullet {\sf b}\equiv   \alpha {\sf b}$, and ${\sf a}\bullet{\sf b}\equiv {\sf a}\cdot {\sf b}$.  Linear extrapolation of the inherent groupoid product (invocation of the distributive law) makes the vector space $\mathbb{R}\oplus \mathbb{R}^{m+n}$ a unital algebra with the multiplicative identity ${\bf 1}$ as simply $1$ (note that, depending on the context, ``1" can denote either a field element or an algebra element).    When $n=0$ and $m\ge 2$, the latter is known as a Spin Factor Jordan Algebra \cite{mccrimmon:2004}.

An element of a unital algebra that has an inverse is a ``unit".  Units in associative algebras and alternative algebras have unique inverses.   Spin Factor Jordan Algebras are neither associative nor alternative, but inverses of units are unique.  
That is, for a Spin Factor Jordan Algebra element $(\alpha +{\sf a})$, suppose it has an inverse $(\beta+{\sf b})$.  Then the expression on the right-hand-side of (\ref{12/24/18.1}) is equal to $1$.  Thus, $\alpha \beta + {\sf a} \cdot {\sf b} =1$ and $ \beta {\sf a} + \alpha {\sf b} = 0$.  Solving these two equations yields the unique inverse $\beta + {\sf b} = \frac{\alpha - {\sf a}}{\alpha^2- {\sf a}\cdot{\sf a}}$.  We define the conjugate operation $^*$ via $(\alpha+{\sf a})^* \equiv (\alpha-{\sf a})$.  Thus, we have the unique inverse of a unit in a Spin Factor Jordan Algebra as,
 \[ (\alpha+{\sf a})^{-1} = \frac{(\alpha+{\sf a})^*}{\alpha^2 - {\sf a}\cdot{\sf a}}.\]

We now proceed to derive the geometry of a Spin Factor Jordan Algebra.
   
   \subsection{The spacetime inner product} In Section 2.1 we derived a particular not associative unital four-dimensional algebra whose units have unique inverses.  Assuming it is even a reasonable task, how does one identify the geometry of such an algebra?  
   
   One can begin with the task of assigning a ``norm" operation on the algebra's elements.  For the case of a real finite-dimensional {\it associative} algebra, each fixed element $x$ can be profitablely viewed as a linear transformation operating on the vector space of algebra elements through the algebra's product operation, i.e., mapping any element $y$ to $xy$.  Accordingly, the ``norm of $x$" can be defined as the determinant of this linear transformation (i.e., the determinant of the element's image under the left regular representation).  This choice is compelling (and standard in ring theory), especially given the geometrical interpretation of the determinant in terms of the volume of an implied parallelpiped.  
    However, the argument in favor of this norm is considerably weakened for the case of {\it not associative} algebras, since the collection of linear transformations associated with each of the algebra elements as above is an associative algebra, while the algebra comprised of the elements themselves is not associative (so there is no algebra homomorphism linking them, unlike the case of an associative algebra and its left regular representation).  Thus, a determinant does not seem to be an appropriate structure to associate with an element of a not associative algebra.  There is a way to salvage the situation for power associative algebras, by taking the norm to be the constant term of the minimal polynomial (which remains well defined), by analogy with the constant term of the characteristic polynomial in the associative algebra case, which is the determinant in the associative algebra setting \cite{jacobson:1963}.  For the case at hand, one can also use the property that a Spin Factor Jordan Algebra can be cast as a $^*$-algebra whose Hermitian elements are real \cite{mccrimmon:2004}.  However, we will pursue a different route, which gets to essentially the same place in this application, but in certain respects is a ``richer" concept than the usual norm on an algebra
    \cite{fgreensite:2023}.
   
   For the remainder of this subsection, ``$\,\cdot\,$" denotes the Euclidean inner product.
      
 We begin by deriving inspiration from a means of identifying a norm on a vector space which has not been given the structure of an algebra.
     For the case of $\mathbb{R}^n$, this entails identifying a function $\ell:\mathbb{R}^n\rightarrow\mathbb{R}$ that encodes features we wish to graft onto $\mathbb{R}^n$, such as that which attends the (degree-1 positive) homogeneity expressed by  $\ell(\alpha s)=\alpha \ell(s)$ for $\alpha >0$, and the concept of a unit sphere $\{s\in\mathbb{R}^n: \ell(s)=1\}$ with respect to which, e.g., isotropy can be defined.  We might also insist that $\ell$ be continuously differentiable on relevant domains, and have an exterior derivative as $\nabla \ell(s)\cdot ds$, where $\nabla \ell(s)$ is defined as the ordered $n$-tuple of coordinate-wise one-dimensional derivatives.  Evidently, $s$ and $\nabla \ell(s)$ are both members of $\mathbb{R}^n$.  However, the latter behave very differently if $s$ is replaced by $\alpha s$.  By the Euler Homogeneous Function Theorem, degree-1 positive homogeneity of $\ell$ implies $s\cdot\nabla \ell(s) = \ell(s)$.  For $\alpha>0$, replacing $s$ by $\alpha s$ then leads to $\alpha s \cdot \nabla \ell(\alpha s) = \ell(\alpha s) = \alpha \ell(s) = \alpha s\cdot \nabla \ell(s)$, implying $\nabla \ell(\alpha s) = \nabla \ell(s)$ (where $\nabla \ell(\alpha s)$ means that $\nabla \ell$ is evaluated at $\alpha s$). On the other hand, for degree-1 positive homogeneous $\ell$, the expressions  $s$ and $\ell(s)\nabla \ell(s)$ {\it do} behave the same way when $s$ is replaced by $\alpha s$ (each expression simply being multiplied by $\alpha$).  So a ``simplest" candidate for $\ell$ presents itself as the solution to, 
  \begin{equation}\label{8/5/22.1} s = \ell(s) \nabla \ell(s),\end{equation} under the constraint,
  \begin{equation}\label{8/5/22.2} \ell(\nabla \ell(s)) = 1,\end{equation} since satisfaction of (\ref{8/5/22.2}) is a necessary condition for a nonnegative solution of (\ref{8/5/22.1}) to be degree-1 positive homogenous (e.g., apply $\ell$ to both sides of (\ref{8/5/22.1}) and apply the homogeneity condition).   It is easily shown that $\ell(s)$ is then the Euclidean norm (derivation of the above two displayed equations from appropriate postulates of Euclidean geometry thereby leads to a novel proof of the Pythagorean Theorem \cite{fgreensite:2022b}).  
  
   But without unduly compromising the simplicity of the rationale, 
  we could propose that (\ref{8/5/22.1}) be replaced by \begin{equation}\label{9/10/22.1}  Ls = \ell(s) \nabla \ell(s) ,\end{equation} for a linear transformation $L:\mathbb{R}^n\rightarrow \mathbb{R}^n$, since we would still have $Ls$ and $\ell(s) \nabla \ell(s)$ behaving the same when $s$ is replaced by $\alpha s$, $\alpha>0$.   Applying the dot product with $s$ to both sides of (\ref{9/10/22.1}), the Euler Homogeneous Function Theorem then implies \begin{equation}\label{9/10/22.2} s'Ls = s \cdot Ls  = \ell^2(s),\end{equation} where $L$ in this equation is understood to be the matrix  associated with the above linear transformation $L$ with respect to the standard basis of $\mathbb{R}^n$.  According to (\ref{9/10/22.1}), $L$ must be a real symmetric matrix since $Ls$ is a gradient (i.e., the right-hand-side of (\ref{9/10/22.1}) is equal to $\frac{1}{2}\nabla\ell^2(s)$). According to (\ref{9/10/22.2}), $L$ must be positive semi-definite.  From the Polarization Identity, it is seen that the above leads to arbitrary inner product spaces on $\mathbb{R}^n$ - a nice, if simple, generalization.
  
  So now we consider the case where vector space $\mathbb{R}^n$ is endowed with the additional structure of a unital algebra and where the inverse of an element is unique if the inverse exists.  In the setting where the algebra is associative (which also guarantees the uniqueness of the inverse of a unit), there is already a compelling rationale for the ``usual norm" as referred to in the second paragraph of this section.  Nevertheless,  one can alternatively speculate how the ``vector space norm treatment" of the prior two paragraphs can incorporate the new opportunity for an element's inverse to influence construction of a ``norm-like" real function $\ell(s)$ on the units of the algebra (something not explicitly addressed by the usual norm).  Proceeding as before, one recognizes that $s^{-1}$ and $\nabla \ell(s)$ are members of $\mathbb{R}^n$, and for $\alpha>0$ we have $\nabla \ell(\alpha s)=\nabla \ell(s)$ if degree-1 positive homogeneity is again mandated.  While the expressions $s^{-1}$ and $\nabla \ell(s)$ thereby behave differently when $\alpha s$ replaces $s$, this time it is $s^{-1}$ and $\frac{\nabla\ell(s)}{\ell(s)}$ that behave the same under that replacement, in that both expressions are simply multiplied by $\frac{1}{\alpha}$ - analogous to the situation in the prior two paragraphs.  So now, based on the treatment applied in the paragraph before last, we might be tempted to equate $s^{-1}$ and $\frac{\nabla\ell(s)}{\ell(s)}$ - except that the latter is the gradient of $\log \ell(s)$, while $s^{-1}$ is in general not a gradient.  That problem resolves if we can find a linear transformation $L$ such that $Ls^{-1}$ \underline{is} a gradient, i.e., satisfies the exterior derivative condition, \begin{equation}\label{12/22/21.1}  d\left([Ls^{-1}]\cdot ds\right) = d\left((ds)'Ls^{-1}\right)=0.\end{equation}  Again, the standard basis on $\mathbb{R}^n$ is assumed, so that $L$ in (\ref{12/22/21.1}) is the matrix associated with the above linear transformation $L$ with respect to the standard basis.

  But there is also the Euler Homogeneous Function Theorem to deal with, whereby $\frac{s\cdot\nabla\ell(s)}{\ell(s)} = 1$. Denoting the multiplicative identity of the algebra as ${\bf 1}$, we can define $\|{\bf 1}\|^2 \equiv \mathbf{1}\cdot\mathbf{1}$.   We are free to constrain $L$ by the requirement that ${\bf 1}\cdot (L{\bf 1}) ={\bf 1}' L{\bf 1} = \|{\bf 1}\|^2$.   It follows that we can now propose, \begin{equation}\label{11/7/20.1} Ls^{-1} = \|{\bf 1}\|^2\frac{\nabla\ell(s)}{\ell(s)} ,\end{equation} if $L$ satisfies both (\ref{12/22/21.1}) and \begin{equation}\label{7/7/21.1} s\cdot Ls^{-1} = s'Ls^{-1} = \|{\bf 1}\|^2,\end{equation} since constraint (\ref{7/7/21.1}) is a necessary condition for a solution of  (\ref{11/7/20.1}) to be degree-1 positive homogeneous (e.g., take the dot product of both sides of  (\ref{11/7/20.1}) with $s$ and apply the Euler Homogeneous Function Theorem).  Since (\ref{7/7/21.1}) must hold if $s$ is replaced by $s^{-1}$ as both are units, we expect $L$ to be a real symmetric matrix.  
  
  Note that (\ref{11/7/20.1}) is analogous to (\ref{9/10/22.1}), and (\ref{7/7/21.1}) is analogous to (\ref{9/10/22.2}).  Indeed, $L$ is akin to an inner product matrix, but a significant feature in this algebra setting compared to the earlier ``naked" vector space setting is that $L$ need not be positive semi-definite.   And,  if $L$ satisfying (\ref{12/22/21.1}) and (\ref{7/7/21.1}) can be found, then one can provide degree-1 positive homogeneous $\ell(s)$ by integrating (\ref{11/7/20.1}), at least in some simply connected neighborhood of ${\bf 1}$.  That is, mandating $\ell(\mathbf{1}) = 1$, we can define \begin{equation}\label{12/19.22.1} \ell(s) \equiv \text{exp}\left(\frac{1}{\|\mathbf{1}\|^2}\int_{\mathbf{1}}^s [Lt^{-1}]\cdot dt\right),\end{equation} where the integral is path-idependent due to (\ref{12/22/21.1}), and real symmetric $L$ must satisfy (\ref{7/7/21.1}) due to the homogeneity requirement.   
  
According to the program of \cite{fgreensite:2023},  $\ell(s)$ is called  a ``unital norm". 
    Real symmetric matrix $L$ satisfying (\ref{12/22/21.1}), but not necessarily (\ref{7/7/21.1}), is called an ``uncurling metric".
  Since $L$ can be considered an inner product matrix, (\ref{7/7/21.1}) can be viewed as a condition mandating that the multiplicative inverses $s,s^{-1}$ in their role as vector space members also act inversely  with respect to the metric $L$.
  
 Lest it be thought that there is something eccentric about this particular identification of a norm with an algebra, we can first note  that for real finite-dimensional unital associative algebras a non-empty family of unital norms always exists and always contains a member ``essentially" equivalent to the usual norm \cite{fgreensite:2023}. 
 
The Spin Factor Jordan Algebra on the spacetime vector space $\mathbb{R}\oplus\mathbb{R}^{3,0}$ introduced in Section 2.1 is not associative, but its units do have unique inverses, and the algebra's multiplicative identity $\mathbf{1}$ is simply 1.  Thus, we can apply the above program to this algebra.  It is easy to show that the only $L$ satisfying (\ref{12/22/21.1}) and (\ref{7/7/21.1}) for this algebra is the conjugate operation, i.e. for $s=\sigma +{\sf s}$, this is $s^* = (\sigma +{\sf s})^* \equiv \sigma -{\sf s}$.  In the nomenclature of \cite{fgreensite:2023}, (\ref{12/19.22.1}) then produces the ``special unital norm", \begin{equation}\label{9/16/22.1} \ell(\sigma+{\sf s}) = \sqrt{\sigma^2-{\sf s}\cdot{\sf s}} = \sqrt{(\sigma+{\sf s})\bullet (\sigma+{\sf s})^*} ,\end{equation}  valid in a some simply connected neighborhood of ${\bf 1}$.  A norm defines an invariant on the space, and so does its square.  The square of the right-hand-side of (\ref{9/16/22.1}) is a quadratic form, which then supplies another invariant on the space by the Polarization Identity, as the inner product,
 \begin{equation}\label{6/29/21.2} \langle \alpha+{\sf a},\beta+{\sf b}\rangle = \frac{(\alpha+{\sf a})\bullet(\beta+{\sf b})^*+ (\beta+{\sf b})\bullet(\alpha+{\sf a})^*}{2} = \alpha\beta-{\sf a}\cdot{\sf b}.\end{equation}  Thus, the Lorentzian geometry of spacetime is derived.
 
 \bigskip
 
 It is worthwhile comparing the above with other algebraic formulations of special relativity, as in the next two subsections. 

\subsection{A non-unital four-dimensional spacetime algebra}
An additional algebra with vector space of elements given by $\mathbb{R}\oplus\mathbb{R}^{3,0}$ can be generated that is closely allied with the Spin Factor Jordan Algebra, but for which the algebra product is the sum of a symmetric bilinear form (given by (\ref{6/29/21.2})) and an antisymmetric bilinear product (in that regard, emulating the geometric product in a Clifford Algebra).  This is accomplished by defining a new algebra product ``$\circ$" on $\mathbb{R}\oplus\mathbb{R}^{3,0}$ as, \begin{eqnarray}\nonumber (\alpha+{\sf a})\circ (\beta+{\sf b}) &\equiv& (\alpha+{\sf a})\bullet (\beta+{\sf b})^* \\ &=&    [\alpha \beta - {\sf a} \cdot {\sf b}] + [ \beta {\sf a} - \alpha {\sf b}] \in \mathbb{R}\oplus \mathbb{R}^{3}.\label{6/29/21.1}\end{eqnarray}  
As demonstrated in Section 3, the algebra defined by (\ref{6/29/21.1}) can be shown to arise in the context of the Geometric Algebra  generated by $\mathbb{R}^{1,3}$.  The relativistic significance of the antisymmetric component of its product is not obvious, since it depends very specifically on the particular $\mathbb{R}^{3,0}$ space of an observer - ultimately deriving from the solipsistic argument of Section 2.1.  However, its broader role will also be clarified in the Geometric Algebra context of Section 3. 

While the Spin Factor Jordan Algebra is a unital algebra, it is easily appreciated that the algebra with product (\ref{6/29/21.1}) is not unital.  The only candidate for the multiplicative identity would be the algebra element $1$.  This is indeed a right identity, but is not a left identity.

\subsection{Comparison with previous formulations of spacetime using Geometric Algebras}
The Geometric Algebra generated by $\mathbb{R}^{1,3}$ is known as the Spacetime Algebra (STA) \cite{hestenes:2003}.  That is, STA is the Clifford Algebra $C\ell_{1,3}(\mathbb{R})$.  The Geometric Algebra generated by $\mathbb{R}^{3,0}$, which is the even subalgebra of $C\ell_{1,3}(\mathbb{R})$, is also of interest, and known as the Algebra of Physical Space (APS) \cite{baylis:1989}. 
APS and STA  each provide an environment suitable for framing relativistic physics.   Although STA assumes Lorentzian 4-space geometry, APS does not (instead, it assumes Euclidean geometry on $\mathbb{R}^3$).  Thus, APS can be thought of as containing an implicit mathematical derivation of the Minkowski inner product on its subspace of paravectors $\mathbb{R}\oplus \mathbb{R}^{3,0}$.  But for that implicit derivation to be possible in the APS context, it must be assumed {\it a priori} that physical units attached to $\mathbb{R}$ (time) and $\mathbb{R}^3$ (physical space) be the same - for otherwise application of an ``isometry"  will in general cause the time components of paravectors to acquire mixed ``time + space" units without evident physical meaning (thereby precluding isometries).   However, a single remaining spacetime physical unit  can still pertain, because of the existence of the Geometric Algebra's multivectors.
     
 In contrast, the derivation of geometry on $\mathbb{R}\oplus \mathbb{R}^{3,0}$ using the Spin Factor Jordan Algebra, according to Sections 2.1 and 2.2, makes the even more radical assumption that the relevant physical units can be replaced by the mathematical unit $1$ (the multiplicative identity of the field $\mathbb{R}$). That is, this program (as presented in Section 2.1) proceeds from the notion that $\mathbb{R}$ and $\mathbb{R}^3$ should be multiplied together (as a tensor product) and added together (as a direct sum) to provide the environment for linearly scaling space and linearly indexing trajectories. The means of identifying geometry on the direct sum is entirely based on generating a self-contained algebra on it - accomplished by linear extrapolation of the groupoid inherent in $\mathbb{R}^{3,0}$.  But the three constituent product operations that are part of the definition of $\mathbb{R}^{3,0}$ could not produce a groupoid if the members of $\mathbb{R}^{3,0}$ are attached to physical units.  This is because the field of scalars  $\mathbb{R}$ must be inherently devoid of physical units (its members only scale), while the product of a member of $\mathbb{R}^{3,0}$ attached to a physical unit with another member of $\mathbb{R}^{3,0}$ attached to a physical unit must map to something attached to the square of the physical unit - and thus cannot map into the field of scalars.  Thus, the groupoid could not exist.
      It follows that physical interpretation and associated physical implications of this program require that the relevant physical unit  be replaceable with the mathematical unit $1$.
    
    However, this prior assertion that the physical units of time and space can be replaced by the mathematical unit is a very sound assumption.  Equivalently, this assumption first predicts a factor equating physical units of time and space - which of course turns out to be the speed of light $c$ - and secondly posits that the remaining physical unit, e.g., meters, ultimately disappears. And indeed, the exchange of the physical unit of length for the mathematical unit ``1" does in fact occur in the gravitational equations (with $\Lambda\ne 0$), as noted in the Introduction.
    
    \bigskip
    
In summary, inclusion of a nonzero Cosmological Constant leads to a presentation of the Einstein equations in which physical units can be interpreted as being superfluous.  This motivates the task of finding a non-physical derivation of spacetime's Lorentzian geometry, as accomplished in Section 2.2.  A bias against non-physical derivations of physical principles is natural, but in our program one is at least left with some novel mathematical notions, as explored in \cite{fgreensite:2023}. 
      Furthermore, it is our view that examination of abstract implications of the elimination of physical units is worthwhile in this case, as these include the ``prediction" of vacuum energy.

    \section{The system of the observer and the observed}
    
    In the approaches of STA, APS, as well as our program in Section 2, the speed of light $c$ has the value 1.  Thus, the relativistic boost equations are    \begin{eqnarray} w'_0 &=& \frac{w_0-v\,w_1}{\sqrt{1-v^2}},\label{12/24/19.1} \\ w'_1 &=& \frac{w_1-v\,w_0}{\sqrt{1-v^2}} \label{12/24/19.2},  \\ w'_2&=& w_2\label{12/24/19.3},  \\ w'_3&=& w_3.   \label{12/24/19.4}  \end{eqnarray}  
     The phenomena of time dilation and boost direction length contraction are immediately recognized, and clearly relate to the $w_0$ and $w_1$ coordinates - but one can't tell which of these two coordinates relate to the time versus the boost directions.  So this is a ``time-boost" symmetry of the system of ``the observer" and ``the observed".  Furthermore, there is a planar subspace of $\mathbb{R}^4$ in which things happen, and a plane orthogonal to it in which nothing happens.    This is, of course, quite distinct from the Lorentz invariance that fundamental equations describing physical forces must satisfy.  In the latter context, there is no prescribed observer-observed pair.  In the former context, the symmetry of the physical equations is broken.
     
     We are going to examine this observer-observed system in the service of exploring the relativistic significance of the antisymmetric portion of the product (\ref{6/29/21.1}).  This will be accomplished in the Clifford Algebra context, i.e., using STA as an environment.
  
  In dealing with STA, one chooses grade-1 elements  $\gamma_\mu$ with $\mu =0,1,2,3$ (i.e., members of $\mathbb{R}^{1,3}$) such that $\gamma_\mu\gamma_\nu = -\gamma_\nu\gamma_\mu$ for $\mu\ne \nu$, where $\gamma_0^2 =1$ and $\gamma_i^2=-1$ for $i=1,2,3$.  These conditions specify that $\{\gamma_\mu\}$ is an orthonormal basis for $\mathbb{R}^{1,3}$, due to the Geometric Algebra feature that the symmetric component of the (geometric) product of grade-1 elements is their inner product in $\mathbb{R}^{1,3}$ - so that the inner product of any two of the particular grade-1 elements above is $ \frac{\gamma_\mu\gamma_\nu+\gamma_\nu\gamma_\mu}{2} = \pm \delta_{\mu\nu}$.    A grade-1 element of STA is notated as, e.g., $a = a^\mu \gamma_\mu$, and we will use the boldface font to represent elements spanned by $\gamma_1,\gamma_2,\gamma_3$, i.e.,  ${\bf a} = a^i \gamma_i$, where the repeated index summation convention is understood in both cases.
  
  The product of grade-1 elements is the so-called geometric product, \begin{equation}\label{12/16/22.1} ab = \langle a,b\rangle + a\wedge b, \end{equation} where $\langle \cdot,\cdot\rangle$ is the Minkowski inner product.  That is,  STA is the Geometric Algebra such that  $\langle a,b\rangle$ is invariant under a Lorentz transformation.  But $a\wedge b$ is {\it not} Lorentz invariant - despite STA being a useful environment in which to frame relativistic physics. This ``anomaly" resolves due to the fact that the antisymmetric wedge product $\wedge$ {\it does} have relativistic significance in expression of the invariant 4-volume, i.e., in the realm of pseudoscalars.
  
In the nomenclature of \cite{hestenes:2003}, the STA element $\gamma_0$ is the ``observer", and one then has the so-called ``spacetime split" determined by the members of  \mbox{$P \equiv \{p\gamma_0: p\in \mathbb{R}^{1,3}\}$}.  As a vector space, $P$ is isomorphic to the inner product space on $\mathbb{R}\oplus\mathbb{R}^{3,0}$ as derived from either the Spin Factor Jordan Algebra via (\ref{12/24/18.1}), or the algebra with product (\ref{6/29/21.1}), since  $\{\gamma_i\gamma_0\}$ is an orthonormal basis for $\mathbb{R}^{3,0}$.  That is, using $\gamma_\mu\gamma_\nu = -\gamma_\nu\gamma_\mu$ for $\mu\ne\nu$, and associativity of the  geometric product, we obtain
 \begin{equation}\label{7/4/21.1} \frac{ (\gamma_i\gamma_0) (\gamma_j\gamma_0) + (\gamma_j\gamma_0) (\gamma_i\gamma_0)}{2} = -\left(\frac{\gamma_i \gamma_j + \gamma_j \gamma_i}{2}\right),\end{equation} where the parenthetical expression on the right-hand-side is the given inner product of basis elements on $\mathbb{R}^{0,3}\subset \mathbb{R}^{1,3}$.  Indeed, $\{\gamma_i\gamma_0\}$ generates the even subalgebra of STA, whose space of grade-1 elements is isomorphic to $\mathbb{R}^{3,0}$ - due to the negative sign attached to the parenthetical (ultimately, inner product) expression on the right-hand-side of (\ref{7/4/21.1}).
 
 As a sum of grade-0 and grade-2 elements of STA, the members of $P$ are subject to the STA product (notated by element juxtaposition).  
But two other relevant products on the vector space of members of $P$ can be introduced.  The first comes from an important observation in STA, which is that for any grade-1 element $a$, associativity of this Geometric Algebra implies $(a\gamma_0)(\gamma_0 a) = aa = \langle a,a\rangle \in \mathbb {R}$.  Since $(a\gamma_0)(\gamma_0 a)$ is a scalar independent of whichever algebra element with unit square is selected as observer $\gamma_0$, it defines the isometries of $P$.  Accordingly, the ``spatial reverse" involution, $(a\gamma_0)^*\equiv \gamma_0 a$, is analogous to the conjugate operation in the Spin Factor Jordan Algebra. For $a$,$b$ grade-1 elements of STA, we thus have a first additional product operation as \begin{equation}\label{6/29/21.4} (a\gamma_0)\star (b\gamma_0) \equiv (a\gamma_0) (b\gamma_0)^* = ab = \langle a,b\rangle + a\wedge b,\end{equation} where the right-hand-side is expression of the (geometric) product of STA grade-1 elements as the sum of the $\mathbb{R}^{1,3}$-inner product and the wedge product.   

The second additional product operation on $P$ comes from the orthogonal projection of the above geometric product to $P$, \begin{eqnarray}\label{3/14/22.1}  (a\gamma_0)\,\circ \,(b\gamma_0) \equiv ab -  {\bf a}\wedge{\bf b} &=& \langle a,b\rangle + [a\wedge b - {\bf a}\wedge{\bf b}] \\ &=& \langle a,b\rangle + \left[b^0{\bf a} - a^0{\bf b}\right]\gamma_0. \label{12/16/23.1}\end{eqnarray}  It is easily seen that ``\,$\circ$\," makes vector space $P$ an algebra.  We can also observe that this  algebra is isomorphic to the algebra on $\mathbb{R}\oplus\mathbb{R}^{3,0}$ defined by (\ref{6/29/21.1}). Either of the bracketed expressions on the right-hand-sides of  (\ref{3/14/22.1}) and (\ref{12/16/23.1}) might be considered to be some kind of a ``partial" wedge product of $a,b$ with respect to a specified observer - and we will think of $\left[b^0{\bf a} - a^0{\bf b}\right]$ in that way.  Along these lines, note that in (\ref{6/29/21.4}) the product $\star$ acting on $P$ is unaffected by an observer $\gamma_0$.  In contrast, product $\circ$ is dependent on $\gamma_0$ - i.e., specific to the space orthogonal to $\gamma_0$ in which the bracketed term on the right-hand-side of (\ref{3/14/22.1}) lives. 

The above self-contained algebra on $P$ can now be used to define a self-contained algebra on $\mathbb{R}^{1,3}$ by simply right-multiplying the right-hand-side of (\ref{12/16/23.1}) by $\gamma_0$ to yield the product, \begin{equation}\label{12/18/23.1}   a\diamond b \equiv \langle a,b\rangle\gamma_0 + \left[b^0{\bf a} - a^0{\bf b}\right].\end{equation}
However, in the context of (\ref{12/24/19.1})-(\ref{12/24/19.4}), there is not only an observer $\gamma_0$ but also an observed $\gamma_1$, though according to those equations we don't know which is which.  That is, for a given 4-vector $a$, we might have \[ a = \underline{a} \gamma_0 +   \bar{a} \gamma_1 + a^2 \gamma_2 + a^3 \gamma_3,\] or instead have, \[ a =  \bar{a} \gamma_0 +   \underline{a}\gamma_1 + a^2 \gamma_2 + a^3 \gamma_3.\]  The analogous situation for a 4-vector $b$ is notated similarly.  This ambiguity leads us to define  two products, \begin{equation}\label{12/16/23.2} a\,\underline{\wedge}\, b \equiv   \left|
\begin{array}{ccc}
\underline{b} & \bar{b}    \\
\underline{a}  &   \bar{a}   \end{array}\right|  \gamma_1 +  \left|
\begin{array}{ccc}
\underline{b} & b^2    \\
\underline{a}  &   a^2   \end{array}\right|  \gamma_2 + \left| 
\begin{array}{ccc}
\underline{b} & b^3    \\
\underline{a}  &   a^3 \end{array}\right|  \gamma_3,\end{equation}
and,
\begin{equation}\label{12/16/23.3} a\,\bar{\wedge}\, b \equiv   \left|
\begin{array}{ccc}
\bar{b}  &  \underline{b}   \\
\bar{a}  & \underline{a}     \end{array}\right|  \gamma_1 +  \left|
\begin{array}{ccc}
\bar{b} & b^2    \\
\bar{a}  &   a^2   \end{array}\right|  \gamma_2 + \left| 
\begin{array}{ccc}
\bar{b} & b^3    \\
\bar{a}  &   a^3 \end{array}\right|  \gamma_3.\end{equation}
Each of the above products is comparable to $\left[b^0{\bf a} - a^0{\bf b}\right]$.  That is, $\left[b^0{\bf a} - a^0{\bf b}\right]$ is the right-hand-side of (\ref{12/16/23.2}) with $\underline{a} = a^0$, \mbox{$\bar{a}= a^1$}, \,\mbox{$\underline{b} = b^0$}, $\bar{b}= b^1$, but is the right-hand-side of (\ref{12/16/23.3}) with $\bar{a} = a^0$, $\underline{a}= a^1$, $\bar{b} = b^0$, $\underline{b}= b^1$.

Like the antisymmetric component $a\wedge b$ of the geometric product (\ref{12/16/22.1}) in the context of $\mathbb{R}^{1,3}$, the antisymmetric product component  $\left[b^0{\bf a} - a^0{\bf b}\right]$ in (\ref{12/18/23.1}) does not have immediately obvious relativistic invariance properties.  However, these exist in the context of the observer-observed system referred to at the beginning of this section, which will be brought out using $a\,\underline{\wedge}\,b$ and $a\,\bar{\wedge}\,b$.  As with $a\wedge b$, recognition of their relativistic role involves the pseudoscalars of STA.

Since the wedge product is distributive and antisymmetric, (\ref{12/16/23.2}), (\ref{12/16/23.3}) imply
\begin{eqnarray}\nonumber 
(a\,\underline{\wedge}\, b)\wedge (a\,\bar{\wedge}\, b) &=&
  \left(\left|
\begin{array}{ccc}
\underline{b} & \bar{b} \\    
\underline{a}  &   \bar{a}   \end{array}\right| \left|
\begin{array}{ccc}
\bar{b} & b^2    \\
\bar{a}  &   a^2   \end{array}\right| - \left|
\begin{array}{ccc}
\bar{b}  &  \underline{b}   \\
\bar{a}  & \underline{a}     \end{array}\right| \left| 
\begin{array}{ccc}
\underline{b} & b^2    \\
\underline{a}  &   a^2 \end{array}\right|   \right)\gamma_1\wedge\gamma_2  
\\ \label{12/16/23.7}  && +
\left(\left|
\begin{array}{ccc}
\underline{b} & \bar{b}    \\
\underline{a}  &   \bar{a}   \end{array}\right|\left| 
\begin{array}{ccc}
\bar{b} & b^3    \\
\bar{a}  &   a^3 \end{array}\right|
  -  \left|
\begin{array}{ccc}
\bar{b}  &  \underline{b}   \\
\bar{a}  & \underline{a}     \end{array}\right|\left|
\begin{array}{ccc}
\underline{b} & b^3    \\
\underline{a}  &   a^3   \end{array}\right|  \right)  \gamma_1\wedge\gamma_3 \\
   && +
\left(\left|
\begin{array}{ccc}
\underline{b}  &  b^2   \\
\underline{a}  &   a^2   \end{array}\right| \left|
\begin{array}{ccc}
\bar{b}  &  b^3   \\
\bar{a}  &   a^3   \end{array}\right|   - \left| 
\begin{array}{ccc}
\bar{b} & b^2    \\
\bar{a}  &   a^2 \end{array}\right|  \left| 
\begin{array}{ccc}
\underline{b} & b^3    \\
\underline{a}  &   a^3 \end{array}\right| \right)\gamma_2\wedge\gamma_3. \nonumber\end{eqnarray}
We then have the volume form/pseudoscalar, 
\begin{multline} \gamma_1\wedge\gamma_2\wedge(a\,\underline{\wedge}\, b)\wedge (a\,\bar{\wedge}\, b)  \\ =
\Bigg( \left| \begin{array}{ccc} 
\underline{b} & b^2   \\
\underline{a}  & a^2  \end{array}\right| \left| \begin{array}{ccc}
\bar{b} & b^3   \\
\bar{a}  & a^3  \end{array}\right|   - \left| \begin{array}{ccc}
\bar{b} & b^2   \\
\bar{a}  & a^2  \end{array}\right|
\left| \begin{array}{ccc}
\underline{b} & b^3   \\
\underline{a}  & a^3  \end{array}\right| \Bigg)\, \gamma_0\wedge\gamma_1\wedge \gamma_2\wedge  \gamma_3  \\
= \left( \left| \begin{array}{ccc}
 \underline{b} & \bar{b}    \\
\underline{a}  &  \bar{a}   \end{array}\right| \gamma_0\wedge\gamma_1 \right) \wedge \left( \left|\begin{array}{ccc}
 b^2 & b^3    \\
a^2  &   a^3   \end{array}\right| \gamma_2\wedge  \gamma_3 \right). \label{12/16/23.5} 
\end{multline}
 The first equality in (\ref{12/16/23.5}) results from (\ref{12/16/23.7}) since $\gamma_1\wedge\gamma_1 = 0$. The second equality in (\ref{12/16/23.5}) arises rather remarkably by multiplying out and then rearranging the determinants in the expression on the left-hand-side of the first equality.  The terms in parentheses on the right-hand-side of the second equality of (\ref{12/16/23.5}) can be thought of as pseudoscalars with respect to the planar subspace defined by the observer and the observed and the planar subspace space orthogonal to it, yielding the respective areas of the parallelograms defined by the orthogonal projections of $a,b$ to those planes - the relevant ``area forms" in the ``observer-observed" context.  A simple calculation reveals that the areas are invariant under application of the (ambiguous) boost equations (\ref{12/24/19.1})-(\ref{12/24/19.4}).

Thus, (\ref{12/16/23.5}) expresses the symmetries of the boost equations (\ref{12/24/19.1})-(\ref{12/24/19.4}).  As we have already noted, these are obviously distinct from the those required to be satisfied by the fundamental equations expressing the physical forces of nature - i.e. Lorentz invariance (the latter equations must be independent of the choice of the observer and independent of the choice of the observed).
That is, the role of the partial wedge product (\ref{12/16/23.2}) (derived from the antisymmetric component of the products (\ref{6/29/21.1}), (\ref{12/16/23.1}), by way of (\ref{12/18/23.1})) has its genesis in the solipsistic argument of Section 2.1, where ``spacetime" in that context is generated from an observer's unique physical space.  While that argument does lead to the Minkowski inner product on the observer's resulting 4-space algebra (Lorentzian geometry), the significance of the antisymmetric product component of the observer's algebra product $\circ$ must be understood in the context of the ambiguity in the implied boost equations (\ref{12/24/19.1})-(\ref{12/24/19.4}).

In other words, the symmetric component of product $\circ$ (the first term on the right-hand-sides of (\ref{6/29/21.1}) and (\ref{12/16/23.1})) is independent of the observer $\gamma_0$ and implies the (Lorentz) invariance referable to the vector space of the algebra's elements, while the antisymmetric component of $\circ$ continues to be tied to the observer and can be analyzed in the context of the boost equations (\ref{12/24/19.1})-(\ref{12/24/19.4}) - the latter equations also being characteristic of special relativity (i.e., when one speaks of an observer and that which is observed).  The relativistic meaning of the antisymmetric product component of $\circ$ emerges in a volume form, as does the relativistic meaning of the antisymmetric product component of the geometric product (\ref{12/16/22.1}) in STA.  However, in contrast to the geometric product, the new product $\circ$ incorporates the spacetime split.

To summarize, the ``observer $\gamma_0$" and the ``observed $\gamma_1$" are  time-like and space-like unit magnitude STA grade-1 elements which specify an observer direction and an observed direction as vectors in $\mathbb{R}^{1,3}$, and thereby also specify the Cartesian plane $\gamma_0\times\gamma_1$, the plane orthogonal to it, and their respective area forms.  However, according to the time-boost symmetry of (\ref{12/24/19.1})-(\ref{12/24/19.4}), for a given 4-vector we don't know which of the first two components is the time versus the boost.  The volume form (\ref{12/16/23.5}) embodies this ambiguity, being (anti)symmetric to an exchange of the observer-observed components.  Furthermore, with respect to an observer, the coefficients of the relevant constituent area forms are invariant under a boost in the direction of the observed  - a reflection of the special relativistic features of time dilation and length contraction.


\begin{thebibliography}{99}

       
          \bibitem{riess:1998}  Riess AG, {\it et al.} (1998) Observational evidence from supernovae for an accelerating universe and a Cosmological Constant. Astronomical Journal 116:1009-1038.
    
    \bibitem{perlmutter:1999}  Perlmutter S {\it et al.} (1999) Measurements of $\Omega$ and $\Lambda$ from 42 high-redshift supernovae. Astrophysical Journal 517:565-586.
    
     
      \bibitem{mccrimmon:2004} McCrimmon K (2004) A taste of Jordan algebras, Springer, New York.
     
              \bibitem{jacobson:1963} Jacobson N  (1963) Generic norm of an algebra. Osaka mathematics Journal, 15(1):25-50.

\bibitem{fgreensite:2023} Greensite F (2023) Novel isomorphism invariants of real algebras, arXiv:2306.14995.          
          
     \bibitem{fgreensite:2022b} Greensite, F (2022), A new proof of the Pythagorean Theorem and generalization of the usual algebra norm. arXiv:2209.14119.
          
         \bibitem{hestenes:2003} Hestenes D (2003) Physics with Geometric Algebra. American Journal of Physics 71:691-714.  

           \bibitem{baylis:1989} Baylis W, Jones G (1989) The Pauli Algebra approach to Special Relativity. Journal of Physics A22:1-16.

 
 \end{thebibliography}
  \end{document}